\documentclass[a4paper,11pt]{amsart}
 
\usepackage{pgf}
\usepackage{amscd,amsmath,amssymb,amstext,amsfonts,amsthm,latexsym}
\usepackage{verbatim}
\usepackage[latin1]{inputenc}
\usepackage{graphicx}
\usepackage{mathrsfs,ulem}
\usepackage[all]{xy}
\usepackage{fancyhdr}
\usepackage{rotating}
\usepackage{url}

 \usepackage[usenames,dvipsnames]{pstricks}
 \usepackage{epsfig}
 \usepackage{pst-grad}

\newtheorem{thm}{Theorem}[section]
\newtheorem*{thm2}{Theorem}

\newtheorem{lem}[thm]{Lemma}

\newtheorem{prop}[thm]{Proposition}

\theoremstyle{remark}
\newtheorem{rem}[thm]{Remark}

\newtheorem*{claim}{Claim}

\theoremstyle{definition}
\newtheorem{defi}[thm]{Definition}

\newcommand{\bbQ}{{\mathbb{Q}}}
\newcommand{\bbZ}{{\mathbb{Z}}}

\newcommand{\bbP}{{\mathbb{P}}}
\newcommand{\bbG}{{\mathbb{G}}}
\newcommand{\bbT}{{\mathbb{T}}}
\newcommand{\bbH}{{\mathbb{H}}}
\newcommand{\bb}{{\mathbb{B}}}
\newcommand{\UU}{\mathcal{U}}

\newcommand{\ud}{\mathrm{d}}

\newcommand{\Id}{\mathrm{Id}}
\newcommand{\Aut}{\mathrm{Aut}}
\newcommand{\Sing}{\mathrm{Sing}}

\newcommand{\ord}{\mathrm{ord}}
\newcommand{\Stab}{\mathrm{Stab}}
\newcommand{\Tors}{\mathrm{Tors}}
\newcommand{\D}{\Delta}

\newcommand{\sts}{\textrm{such that }} 
\newcommand{\qe}{\textrm{q.e.}\,\,}
\newcommand{\ie}{\textrm{i.e.}\,\,}
\newcommand{\cf}{\textrm{cf.}\,\,}

\newcommand{\T}{\tilde \tau}

\newcommand{\no}[1]{\textit{#1}}
\newcommand{\ov}[1]{\overline{#1}}
\newcommand{\q}[1]{``{#1}''}
\newcommand{\tr}[1]{\langle{#1}\rangle}

\title[Mixed quasi-\'etale surfaces with $p_g=0$]{Mixed quasi-\'etale surfaces,
 new surfaces of general type with $p_g=0$ and their fundamental group} 																											
\makeatletter

\@addtoreset{equation}{section}
\makeatother
\author[D. Frapporti]{DAVIDE FRAPPORTI}
\keywords{Surfaces of general type, finite group actions} 
\subjclass[2000]{14J29, 	14Q10 , 14Q99, 20F34, 20F05, 58E40} 

\date{\today}

\begin{document}

\begin{abstract} We call a projective surface $X$ \no{mixed quasi-\'etale quotient}
if there exists a curve $C$ of genus $g(C)\geq 2$ and a finite group $G$ that acts on $C\times C$ exchanging the factors 
\sts  $X=(C\times C)/G$ and the map $C\times C \rightarrow X$ has finite branch locus.
The minimal resolution of its singularities is called \no{mixed quasi-\'etale surface}. 
We study the mixed quasi-\'etale surfaces under the assumption that
$(C\times C)/G^0$ has only nodes as singularities, where $G^0\triangleleft G$ is the index two subgroup
 of the elements that do not exchange the factors. 

We classify the minimal regular surfaces with $p_g=0$ whose canonical model is a 
mixed quasi-\'etale quotient as above. All these surfaces are of general type and 
as an important byproduct, we provide an example of a numerical Campedelli surface with
topological fundamental group $\bbZ_4$, and we realize 2 new topological types
of surfaces of general type. Three of the families we construct are $\bbQ$-homology projective planes. 
\end{abstract}
\maketitle

\section*{Introduction}

It is a well known fact that each Riemann surface with $p_g=0$ is isomorphic to $\bbP^1$.
At the end of XIX century M. Noether conjectured that an analogous statement holds for 
the surfaces: in modern words he conjectured that every  smooth projective surface with $p_g=q=0$ 
be rational.
The first counterexample to this conjecture is due to F. Enriques (1896),
he introduced the so called Enriques surfaces (see \cite{Enriques}), 
that are surfaces of special type.
The first examples of surfaces of general type with $p_g=0$ have been constructed in the 30's 
by L. Campedelli and  L. Godeaux.

The idea of Godeaux to construct surfaces was to consider the quotient
of simpler surfaces by the free action of a finite group.
In this spirit, Beauville (\cite{Beau}) proposed a simple construction
of a surface of general type, considering the 
quotient $(C\times C)/G$ where $C$ is the Fermat plane quintic and 
$G$ is the finite group $(\bbZ_5)^2$ that acts freely on the product.
This construction leads to a surface with $p_g=q=0$ and $K^2=8$.

Nowadays, some example of surfaces of general type with $p_g=0$ are known (see \cite{Survey} for a more
detailed discussion), but the classification of them is far from being complete. 
Generalizing the Beauville example, we can consider the quotient $(C_1\times C_2)/G$, 
where the $C_i$'s are  Riemann surfaces of genus at least two, and $G$ is a finite group.
By \cite{Cat00}, there are two cases: the \no{mixed case} where  the 
action of $G$ exchanges the two factors (and then $C_1\cong C_2$); 
and the \no{unmixed case} where $G$ acts diagonally.

After \cite{Cat00} many authors started studying the surfaces
that appear as quotient of a product of curve, see
\cite{BC04}, \cite{BCG08}, \cite{BCGP08} and \cite{BP10} for $p_g=q=0$; 
\cite{CP09}, \cite{Pol07},\cite{Pol09} and \cite{MP10} for $p_g=q=1$;
\cite{Penegini} for $p_g=q=2$.  
In all these articles the authors work either in the unmixed case or
in the mixed case under the assumption that the group acts freely.

The main purpose of this article is to extend the results and the strategies
of the above mentioned cases in the non free mixed case.\\
Let $C$ be a Riemann surface of genus $g\geq 2$, let $G$ be a finite group that
acts on $C\times C$ with a mixed action, and let $G^0\triangleleft G$ be the index two subgroup of the
elements that do not exchange the factors. We say that $X=(C\times C)/G$ is a \no{mixed quasi-\'etale quotient} if the 
quotient map $C\times C \rightarrow (C\times C)/G$ has finite branch locus (see Definition \ref{MIX});  
the minimal resolution of its singularities is a \no{mixed quasi-\'etale surface} .
In this paper we assume that $(C\times C)/G^0$ has only nodes as singularities and we will construct 
some new surfaces of general type with $p_g=0$ as 
mixed quasi-\'etale surfaces. In particular we prove the following:

\begin{thm2}
Let $S$ be a minimal regular surface with $p_g(S)=0$ whose canonical
model is the mixed \qe quotient
 $X=(C\times C)/G$ \sts $(C\times C)/G^0$ has at most nodes as singularities.
Then $S$ is of general type and belongs to one of the 13 families collected in Table \ref{tab1}.
\end{thm2}


\small{
\begin{table*}[!h]
\centering
\hspace{-2.5cm}
\begin{minipage}{\textwidth}

\begin{tabular}{|c|c|c|c|c|c|c|c|c|}
\hline

$K^2_S$\hspace{-.1cm} & $\Sing(X)$ & Type & $G^0$ & $G$& $b_2$  & $H_1(S,\bbZ)$ & $\pi_1(S)$ & \hspace{-.1cm} Label\hspace{-.1cm}\\

\hline\hline

1	& $2 \, A_1,2\,  A_3$ & $2^3,4$ & $D_4\times  \bbZ_2$ & $\bbZ_2^3\rtimes \bbZ_4$ & $1$& $\bbZ_4$ & $\bbZ_4$  &  {A.1.1}  \\

\hline\hline

2 & $6 \, A_1$ & $2^5$ & $\bbZ_2^3$ & $\bbZ_2^2\rtimes \bbZ_4$& $2$& $\bbZ_2\times\bbZ_4$ & $\bbZ_2\times\bbZ_4$ &{A.2.1} \\

2 & $6 \, A_1$ & $4^3$ & $(\bbZ_2\times \bbZ_4)\rtimes\bbZ_4$ & G(64, 82)& $2$ & $\bbZ_2^3$&$\bbZ_2^3$&{A.2.2}\\

2 & $  A_1,2\, A_3$ & $2^3,4$ & $\bbZ_2^4\rtimes \bbZ_2$ & $\bbZ_2^4\rtimes \bbZ_4$& $1$  & $\bbZ_4$ & $\bbZ_4$&{A.3.1}\\

2 & $  A_1,2\,  A_3$ & $2^2,3^2$ & $\bbZ_3^2\rtimes \bbZ_2$& $\bbZ_3^2\rtimes \bbZ_4$& $1$ & $\bbZ_3$ &$\bbZ_3$& {A.3.2}\\

\hline\hline

4 & $4 \, A_1$ & $2^5$ & $D_4 \times \bbZ_2$& $D_{2,8,5}\rtimes \bbZ_2$ & $2$& $\bbZ_2\times \bbZ_8$ &$\bbZ_2^2\rtimes \bbZ_8$&{A.4.1}  \\

4 & $4 \, A_1$  & $2^5$ & $\bbZ_2^4$ & $(\bbZ_2^2\rtimes \bbZ_4) \times \bbZ_2$ & $2$  & $\bbZ_2^3\times \bbZ_4$& K-N &{A.4.2}   \\

4 & $4 \, A_1$  & $4^3$ &  G(64, 23) & G(128, 836)& $2$ & $\bbZ_2^3$& $\bbZ_4^2\rtimes \bbZ_2$& {A.4.3} \\

\hline\hline

8 & $\emptyset$ & $2^5$ & $D_4 \times \bbZ_2^2$ &$(D_{2,8,5}\rtimes \bbZ_2) \times \bbZ_2$& $2$  &$\bbZ_2^3\times \bbZ_8$
& $\infty$
&{A.5.1} \\

8 & $\emptyset$ & $4^3$ & G(128, 36) & G(256, 3678)& $2$ & $\bbZ_4^3$  
& $\infty$
& {A.5.2}  \\

8 & $\emptyset$ & $4^3$ & G(128, 36) & G(256, 3678)& $2$ & $\bbZ_2^4\times \bbZ_4$  
& $\infty$
& {A.5.3} \\

8 & $\emptyset$ & $4^3$ & G(128, 36) & G(256, 3678)& $2$ &$\bbZ_2^2\times \bbZ_4^2$  
& $\infty$
 & {A.5.4} \\

8 & $\emptyset$ & $4^3$ & G(128, 36) & G(256, 3679)& $2$ & $\bbZ_2^2\times \bbZ_4^2$  
& $\infty$
& {A.5.5}  \\
\hline
\end{tabular}
 
\caption{The surfaces} \label{tab1}
\end{minipage}
\end{table*}
}

In the first column of the table we report the $K_S^2$ of the surface; $\Sing(X)$ represents the 
singularities of $X$.
The column Type gives the type of the set of spherical generators (see Definition \ref{ssg}) in a compacted way, e.g. $2^3,4=(2,2,2,4)$.
The columns $G$ and $G^0$ give the group and its index two subgroup.
The column $b_2$ reports the second Betti number of the canonical model $X$.
The last two columns report the first homology group and the topological fundamental group of the surface.

For the groups occurring in the paper we use the following notation:
we denote by $\bbZ_d$ the cyclic group with $d$ elements.
$D_{p,q,r}$ is the generalized dihedral group with presentation:
$D_{p,q,r} = \tr{x, y|x^p, y^q, xyx^{-1}y^{-r}}$ and $D_n := D_{2,n,-1}$ is the dihedral group of order $2n$.
$G(a, b)$ denotes  the $b$-th group of order $a$ in the MAGMA database of groups.
An expanded version of Table \ref{tab1}, can be downloaded from:
$$\mbox{\url{http://www.science.unitn.it/~frapporti/papers/surfaces1.pdf}}\,$$

In the $K^2=8$ case, the quotient map $C\times C \rightarrow (C\times C)/G$ is \'etale (i.e. $G$ acts freely),
this case has been already classified by \cite{BCG08}, see also Remark \ref{missing}. 
 
We point out that in Table \ref{tab1} appears a numerical Campedelli surface ($K^2=2$) 
with  topological fundamental group (and therefore algebraic fundamental group) 	$\bbZ_4$. 
By the works of M. Reid and others (\cite{Rei}, \cite{MP08}, \cite{MPR09})
it is known that the algebraic fundamental group of a numerical Campedelli surface 
is either abelian of order $\leq 9$ or it is the quaternion group $Q_8$.
The question whether all these groups occur has been open for a while, and 
a similar question for the topological fundamental group has been posed in \cite[Question 2.17]{Survey}.
The answer to the question for the algebraic fundamental group is affirmative. 
Indeed, the last open case, $\bbZ_4$, is realized by our example and by a completely
different construction found independently by \cite{PPS10}.
We note that the topological fundamental group of \cite{PPS10} is not known.
We mention (cf. \cite{Survey}) that after our construction, the only open case
left for the question on the topological fundamental group is $\bbZ_6$.

We note that our constructions provide at least other 2 topological types of surfaces 
which were not known before, see Remark \ref{news}.

We also note that three of the families we construct are 
$\bbQ$-homology projective planes in sense of \cite{HK} and \cite{Keum10}, see Remark \ref{qhpp}.

\

The paper is organized as follows: in Section \ref{Notation} we set some notation. 
In Section \ref{mixed} we introduce the mixed quasi-\'etale surfaces and we
investigate them and their singularities under the above assumptions.
In Section \ref{Costruzione} 
we relate the numerical invariants $e$ and $K^2$ with the singularities of $X$, the order of $G$ and the genus
of $C$.
In Section \ref{Pi1} we explain how to calculate the
fundamental group of the surfaces using Armstrong's results (\cite{Arm65}, \cite{Arm68}).
Section \ref{class} is devoted to the proof of the main theorem of the paper.
We also compare our surfaces with the constructions in literature.

\section{Notation}\label{Notation}

We will use the same notation as in \cite{BCGP08}.
Given natural numbers $m_1, \ldots, m_r > 1$ the 
\no{polygonal group of signature $(m_1, \ldots, m_r)$} is defined as:
\begin{equation}\label{poly}
\bbT(m_1,\ldots ,m_r):=
\tr{c_1, \ldots, c_r \mid c_1^{m_1}, \ldots, c_r^{m_r}, c_1 \cdots c_r}\,.
\end{equation}
Let $H$ be a finite group, we say that an homomorphism $$\psi\colon\bbT(m_1,\ldots ,m_r)\rightarrow H$$ is 
an \no{appropriate orbifold homomorphism} if it is surjective and $h_i:=\psi(c_i)$ has order $m_i$.

\begin{defi}\label{ssg}
Let $H$ be a finite group. A \no{spherical system of generators} of $H$ of type (or signature) $(m_1,\ldots,m_r)$
is a  set of generators $\{h_1,\ldots, h_r\}$ of $H$ \sts $h_1\cdots h_r=1$ and there 
exists a permutation $\sigma \in \mathfrak{S}_r$ such that $\ord(h_i)=m_{\sigma(i)}$ for $i=1,\ldots, r$.
\end{defi}

By Riemann's existence theorem (see \cite{Survey}), any curve $C$ together with an action of
a finite group $H$ on it such that $C/H\cong \bbP^1$ is determined (modulo automorphisms)
by the following data:
\begin{enumerate}
	\item the branch point set $\{p_1, \ldots, p_r\}\subset \bbP^1$;
	\item $\gamma_1,\ldots,\gamma_r\in \pi_1(\bbP^1\setminus\{p_1, \ldots, p_r\})$, where
	 each $\gamma_i$ is a simple geometric loop around $p_i$ and $\gamma_1\cdot\ldots\cdot\gamma_r=1$
	\item an appropriate orbifold homomorphism $\psi\colon\bbT(m_1,\ldots ,m_r)\rightarrow H$	
	with the property that \no{Hurwitz's formula} holds:
	\begin{equation}
	2g-2=|H|\bigg(-2+\sum_{i=1}^r\frac{m_i-1}{m_i}\bigg)\,.
	\end{equation}
	\end{enumerate}

\section{On mixed quasi-\'etale surfaces}\label{mixed}

We start defining the objects of our analysis:

\begin{defi}[\cf {\cite[Proposition 3.15]{Cat00}}]
Let $C$ be a Riemann surface of genus $g(C)\geq 2$, and let $G$ be a finite group.
A \no{mixed action} of $G$ on $C\times C$ is a monomorphism $G\hookrightarrow\Aut(C\times C)\cong \Aut(C)^2\rtimes \bbZ_2$
whose image is not contained in $\Aut(C)^2$.
Given a mixed action we will denote by $G^0\triangleleft G$ the index two subgroup $G \cap \Aut(C)^2$.
A mixed action is \no{minimal} if $G^0$ acts faithfully on both factors.
\end{defi}

\begin{defi}
A \no{mixed quotient} is a surface which arises as quotient $X:=(C\times C)/G$
by a mixed action of $G$ on $C\times C$.
\end{defi}

\begin{rem}[{\cf \cite[Remark 3.10, Proposition 3.13]{Cat00}}]
Every mixed quotient $X$ may be obtained by an unique minimal mixed action.

Let $K_2\times \Id:=G^0\cap(\Aut(C)\times \Id)$ and $\Id\times K_1:=G^0\cap(\Id\times \Aut(C))$, then
$K_1\cong K_2$ as subgroups of $\Aut(C)$.
Moreover $K_1\times K_1$ is a normal subgroup of $G$, and $G/(K_1\times K_1)$ acts mixed and minimally on $(C/K_1)\times (C/K_1)$.
The proof of the uniqueness is analogous to the proof of \cite[Proposition 3.13]{Cat00}.
\end{rem}
	
\begin{defi}
Let $X$ be a mixed quotient. By the previous remark we may obtain $X$ as $C\times C/G$ 
by a minimal mixed action; we will call the map $C\times C\rightarrow X$ \no{the quotient map} of $X$.  
\end{defi}

\begin{rem}
The quotient map can be factorized as follows:
$$C\times C \stackrel{\sigma}{\longrightarrow} Y:=(C\times C)/G^0\stackrel{\pi}{\longrightarrow}X\,.$$
\end{rem}

F. Catanese in \cite[Proposition 3.16]{Cat00} 
gives the following description of minimal mixed actions:

\begin{thm}[{\cite[Proposition 3.16]{Cat00}}]\label{thmix}

Let $G\subseteq \Aut(C\times C)$ be a minimal mixed action. 
Fix $\tau' \in G\setminus G^0$; it determines an element $\tau:=\tau'^2 \in G^0$
and an element $\varphi\in \Aut(G^0)$ defined by $\varphi(h):=\tau' h \tau'^{-1}$.
Then, up to a coordinate change, $G$ acts as follows:
\begin{equation}\label{action}
\begin{split}
g(x, y) &= (gx, \varphi(g)y)\\ 
\tau'g(x, y) &=(\varphi(g)y, \tau g \,x)
\end{split}\qquad for \,\,g \in G^0
\end{equation}

Conversely, for every $G^0\subseteq\Aut(C)$ and $G$ extension of degree 2 of $G^0$, fixed 
$\tau'\in G\setminus G^0$ and defined $\tau$ and $\varphi$ as above, (\ref{action}) defines a minimal
mixed action on $C\times C$.
\end{thm}

F. Catanese gives also a characterization of the mixed quotient whose 
quotient map is \'etale. 
In the following we generalize that statement to the case when the quotient map is quasi-\'etale (see \cite{Cat07}),
 \ie the branch locus is finite.

\begin{thm}\label{nonsplit}

Let $X$ be a mixed quotient provided by a minimal 
mixed action of $G$ on $C\times C$.
The quotient map $C\times C\rightarrow X$ is quasi-\'etale if and only if the exact sequence
\begin{equation}\label{ext}
1 \longrightarrow G^0 \longrightarrow G \longrightarrow \bbZ_2 \longrightarrow 1 
\end{equation}
does not split.

Moreover, if the quotient map is quasi-\'etale, then $\Sing(X)= \pi(\Sing(Y))$.
				
\end{thm}

\begin{proof}

$(\Rightarrow)$
	If there exists $h \in G^0$ 
	\sts $(\tau' h)^2=1$, i.e. $\varphi(h)\tau h=1$, then we get $$\tau' h(x, \tau h x)=(\varphi(h)\tau h x, \tau h x)=(x, \tau h x)\,,$$
	hence the quotient map $C\times C\rightarrow X$ is  ramified along a curve.

$(\Leftarrow)$
We factor the quotient map of $X:=(C\times C)/G$ as 
$$C\times C \stackrel{\sigma}{\longrightarrow} Y:=(C\times C)/G^0\stackrel{\pi}{\longrightarrow}X\,.$$
Since $G^0$ acts faithfully, $\sigma$ is branched only in a finite number of points: $r_1,\ldots, r_t$.
 	Aiming for a contradiction we assume that 
 	there exists a curve $D \subseteq X$ \sts $|\pi^{-1}(q)|=1$ for all $q\in D$.\\ 		
 	Let $q \in D$ \sts $\pi^{-1}(q)=p'\not \in \{r_1,\ldots, r_t\}$. 
 	Since $\sigma$ is a $|G^0|=:n$ to 1 map, we get $\sigma^{-1}(p')=\{p_1, \ldots , p_{n}\}$
 	and $|(\pi \circ \sigma)^{-1}(q)|=n$. It follows that
 	$|\Stab(p_1)|=2$, hence $\Stab(p_1)\cong\bbZ_2$ is generated by an element not in $G^0$.
 	Then (\ref{ext}) splits, a contradiction.

\

Let $\{r_1,\ldots, r_t\}$ be the singular locus of $Y$. If $q\in \Sing(X)\setminus \pi(\Sing(Y))$ then 
 $\pi^{-1}(q)=p'\not \in \{r_1,\ldots, r_t\}$ and we can argue as before to get a contradiction. 
 Then $\Sing(X)\subseteq\pi(\Sing(Y))$, the opposite inclusion is a special case of \cite[Remark 3.1]{Cat07}.
\end{proof}

\begin{defi}\label{MIX}
A \no{mixed quasi-\'etale quotient} $X=(C\times C)/G$ is a mixed quotient with quotient map quasi-\'etale and provided by the corresponding minimal mixed action, as described in Theorem \ref{thmix}.
The minimal resolution of its singularities is called \no{mixed quasi-\'etale surface}.
\end{defi}

\begin{lem}
Let $S\rightarrow X=(C\times C)/G $ be a mixed \qe surface.
Then $q(S)$  equals the genus of $C':=C/G^0$.
\end{lem}
\begin{proof}
From \cite[Proposition 3.5]{MP10} (see also \cite{Freitag}), we have that $$H^0(\Omega_S^1)=(H^0(\Omega_{C\times C}^1))^G\,.$$
Arguing as in \cite[Proposition 3.15]{Cat00}:
$$\begin{array}{lcl}
H^0(\Omega_S^1)&=&(H^0(\Omega_{C\times C}^1))^G =(H^0(\Omega_C^1)\oplus H^0(\Omega_C^1))^G		\\
							 &=&(H^0(\Omega_C^1)^{G^0}\oplus H^0(\Omega_C^1)^{G^0})^{G/{G^0}} \\
							 &=&(H^0(\Omega_{C'}^1)\oplus H^0(\Omega_{C'}^1))^{G/{G^0}}.
\end{array}
$$
Since $X$ is a mixed quotient, $\bbZ_2=G/{G^0}$ exchange the last summands,
hence $q(S)=h^0(\Omega_S^1)=h^0(\Omega_{C'}^1)=g(C')$.
\end{proof}

In \cite{BCG08} the authors have constructed surfaces of general type with $p_g=0$ as mixed quotient with
 \'etale quotient map. We want to use Theorem \ref{nonsplit} to extend their construction, so we will assume 
that the quotient map has finite branch locus.
We further assume that $Y=(C\times C)/G^0$ has only nodes (Du Val singularities of type $A_1$) as singularities.

In a forthcoming paper we will drop this assumption and we will assume that $Y$ has arbitrary singularities.

\begin{prop} \label{sing}
Let $X=(C\times C)/G$ be a mixed \qe quotient. 
Let $p\in \Sing(Y)$ be a singularity of type $A_1$. Then
 $\pi(p)$ is a point of type:
$A_1$ if $p$ is not a ramification point of $\pi$; $A_3$ otherwise.
 \end{prop}

\begin{proof}
Let $(x,y)\in C\times C$ \sts $p=\sigma(x,y)$ is a node in $Y$;
if $p$ is a ramification point of $\pi$, $p$ is fixed by the involution
induced by $G$ on $Y$.
By \cite[Theorem 2.2]{Cat87} and \cite[Theorem 2.5]{Cat87}
the quotient of a node by an involution with isolated fixed points is either 
a point of type $A_3$ or a singular point of type $\frac{1}{4}(1,1)$ ($Y_1$ in the notation of \cite{Cat87}).
We show that the latter case does not happen. 
Let $\eta$  be a generator of $\Stab_G(x,y)\cong \bbZ_4$.
Then $\ud \eta_{(x,y)}= i \cdot Id$ has trace $2i$.
On the other hand, $\eta\notin G^{0}$ so it exchanges the two factors
and then $\ud \eta_{(x,y)}=\left(
		\begin{array}{cc}
		0 & b\\
		c &0
		\end{array}
		\right)$;
in particular it has trace $0$, a contradiction.				
\end{proof}

\section{Constructing surfaces}\label{Costruzione}

\begin{rem}\label{algdata}
A regular mixed \qe surface is completely determined by the following algebraic data:
\begin{itemize}
  \item some points $\{p_1,\ldots, p_r\}\subset \bbP^1$ and
  $\gamma_1,\ldots,\gamma_r\in \pi_1(\bbP^1\setminus\{p_1, \ldots, p_r\})$, where
	 each $\gamma_i$ is a simple geometric loop around $p_i$ and $\gamma_1\cdot\ldots\cdot\gamma_r=1$
	\item a finite group $G$;
	\item a spherical system of generators $(h_1,\ldots,h_r)$ of type $(m_1, \ldots, m_r)$ of an index two subgroup
	$G^0\triangleleft G$ \sts $1 \rightarrow G^0 \rightarrow G \rightarrow \bbZ_2 \rightarrow 1 $
	does not split.
\end{itemize}

\noindent Once we fix $\{p_1,\ldots, p_r\}$, $\{\gamma_1,\ldots,\gamma_r\}$, $G$ and  $(h_1,\ldots,h_r)$, by Riemann's existence theorem we get
a curve $C$ \sts the cover $c\colon C \rightarrow C/G^0\cong \bbP^1$
is branched over $\{p_1, \ldots, p_r \}\subseteq \bbP^1$. 
Using Theorem \ref{thmix} we define a mixed action on $C\times C$ and 
by Theorem \ref{nonsplit} the quotient map is quasi-\'etale.

We note that a mixed \qe surface is determined up to the choice of $r$ points in $\bbP^1$,
hence we get a family of surfaces parametrized by $r-3$ parameters.
\end{rem}

\begin{rem}\label{cl}
Different algebraic data may determine deformation equivalent surfaces. 
Let $\bb_r$ the braid group on $r$ elements and 
consider the action of $\bb_r\times \Aut(G)$ 
on the sets of spherical generators of length $r$:
\begin{equation}\label{hm}
(\gamma,\eta)\cdot(G^0,T):=(\eta(G^0),\eta(\gamma(T)))\,.
\end{equation}
This group action was introduced in 
\cite[Section 1-2]{BCG08} 
(see also \cite[Section 5.1-5.2]{CSGT}), where it is shown that
two pairs $(G^0,T)$ in the same orbit (fixed the branch points) give 
deformation equivalent mixed \qe surfaces. 
\end{rem}

\begin{defi}
Let $X=(C\times C)/G$ be a mixed \qe quotient. If $Y$ is nodal by Proposition \ref{sing} all singularities of $X$ are either of type $A_1$ or of type $A_3$.
We will denote by $s$ the number of nodes and by $t$ the number of $A_3$ singularities of $X$.
\end{defi}

\begin{lem}\label{K2value}
Let $S\rightarrow X=(C\times C)/G$ be a mixed \qe surface.
If $Y$ is nodal then $$K^2_S=K^2_X=8\chi(S)-s-\dfrac{5}{2}t>0\,.$$
In particular $t$ is even.
\end{lem}

\begin{proof}
Arguing as in \cite{BCGP08}, we get 
$$K^2_S=K^2_X=\dfrac{8(g-1)^2}{|G|}>0$$
since $X$ has only canonical singularities and the quotient map
$\pi\circ \sigma$ is quasi-\'etale.
We also get 
$$e(S)=e(X)+s+3t=\dfrac{4(g-1)^2}{|G|}+\dfrac{3s}{2}+\dfrac{15t}{4}\,.$$

\noindent By Noether's formula:
$$12\chi(S)=K^2_S+e(S)=\underbrace{\dfrac{12(g-1)^2}{|G|}}_{3 K_S^2/2}+\dfrac{3s}{2}+\dfrac{15t}{4} \Longrightarrow K_S^2=8\chi(S)-s-\frac{5t}{2}\,.$$
\end{proof}

 We use the combinatorial restriction forced by the assumptions
in o	rder to determine the possible signature $(m_1,\ldots, m_r)$ of the polygonal group, 
and the possible cardinalities of $G$.
We start defining the following numbers:
$$\Theta:=-2+ \sum_{i=1}^r\dfrac{m_i-1}{m_i}\,, \qquad \beta:=\dfrac{K^2_S}{2\Theta}\,.$$
		
\begin{prop} \label{orderG}
Let $S\rightarrow X=(C\times C)/G$ be a regular mixed \qe surface.
Let $(m_1,\ldots, m_r)$ be the signature of the spherical system of generators of $G^0$ associated to $X$; then
$$\Theta >0\,, \qquad\beta=g(C)-1\, \mbox{ and }\quad |G^0|=\dfrac{4\beta^2}{K_S^2}\,.$$
 
\noindent Moreover, each $m_i$ divides $2\beta$ and there are at most $\frac{n}{2}$ indices $i\in \{1, \ldots, r\}$ such that 
$m_i$ does not divide $\beta$, where $n$ is the number of nodes on $Y$.
 \end{prop}
 
\begin{proof}
Let $g$ be the genus of $C$. Since $C/G^0\cong \bbP^1$, by Hurwitz's formula we get 
$$2(g-1)=|G^0|\cdot \Theta \,,$$
hence $\Theta=\dfrac{2(g-1)}{|G^0|}>0$, since $g\geq 2$.
 $K_S^2=\dfrac{4(g-1)^2}{|G^0|}$, therefore 
$$\beta=\dfrac{4(g-1)^2}{2\Theta\cdot |G^0|}=g-1
\Longrightarrow|G^0|=\dfrac{4\beta^2}{K^2_S}\,.$$

\noindent The last claims follow from \cite[Lemma 5.8]{BCGP08}.
\end{proof}

\begin{lem}\label{inequalities}
Under the same assumption of Proposition \ref{orderG}, the following hold:
\begin{itemize}
	\item[a)] $r< K_S^2+5$;
	\item[b)] for all $i$, we have that $m_i \leq \dfrac{1}{M}(K_S^2+1)$, where $M:=\max\{\frac{1}{6},\frac{r-3}{2}\}$.
\end{itemize}
\end{lem}

\begin{proof}

$a)$ Assume that $r\geq K_S^2+5>5$, hence $\Theta\geq -2+\frac{r}{2}=\frac{r-4}{2}>0$.
We get: $$1 \leq \beta=\dfrac{K_S^2}{2 \Theta}\leq \dfrac{K^2_S}{2}\cdot
\dfrac{2}{r-4}=\frac{K_S^2}{r-4}\leq \dfrac{K^2_S}{K_S^2+1}<1\,.$$

$b)$ We can assume $m_1\geq m_i$ for all $i$. Since $\Theta$ is strictly positive then $r\geq 3$.
If $r=3$ at most one $m_i$ can be equal to $2$, hence
$$\Theta+\frac{1}{m_1}=1-\dfrac{1}{m_2}-\dfrac{1}{m_3}\geq \frac{1}{6}\,.$$

\noindent If $r>3$, since $\displaystyle{\Theta=(r-2)-\sum_{i=1}^r\frac{1}{m_i}}$, 
it holds:
$$
\Theta+\frac{1}{m_1}=(r-2)-\sum_{i=2}^r\frac{1}{m_i} 
\geq (r-2)-\frac{r-1}{2}=\frac{r-3}{2}\,.
$$

\noindent Hence $\Theta+\frac{1}{m_1}\geq\max\{\frac{1}{6},\frac{r-3}{2}\}=M$. Since $m_i\leq 2\beta$ 
$$m_1\leq \frac{1}{M}(\Theta \cdot m_1+1)\leq \frac{1}{M}(\Theta \cdot 2\beta +1)=\frac{1}{M}(K_S^2+1)\,.$$
\end{proof}

\subsection{How to count the singularities}

In order to implement an algorithm to construct regular mixed \qe surfaces, 
we need to understand how to count the singularities 
of $Y$ and $X$ starting from the algebraic data (see Remark \ref{algdata}).

\begin{rem}
We recall that the points in $c^{-1}(p_i)$ are the only ones with non trivial stabilizer with respect
to the action of $G^0$ on $C$ and they are
in bijection with the left cosets $\{g K_i\}$, where $K_i:=\tr{h_i}$. 
Note that the point $g K_i$ has stabilizer $gK_i g^{-1}$ and that $|c^{-1}(p_i)|=\dfrac{|G^0|}{m_i}$.
Let $Q\colon Y\rightarrow\bbP^1\times\bbP^1$ be  the map $Q(\sigma(x,y))=(c(x),c(y))$.
\end{rem}

\begin{prop}
Let $X$ be a mixed \qe quotient determined by $(h_1,\ldots,h_r)$ and $G$.
Assume that $Y=(C\times C)/G^0$ has only nodes as singularities.\\
If $m_i=o(h_i)$ is even, let  
	  $d_i:=m_i/2$,
	   $e_i:=|\{g h_i^{d_i}g^{-1}\}_{g \in G^0}|$,
   $$\delta_{ij}:=\left\{ \begin{array}{ll}
1 & \mbox{ if } h_i^{d_i} \mbox{ is conjugated (in $G^0)$ to } \varphi^{-1}(h_j^{d_j})\\
0 & \mbox{ otherwise}
\end{array}
\right.$$
  and $F_i:=\{\tau'\eta \in G\setminus G^0 \mid (\tau'\eta)^2\mbox{ is conjugate in $G^0$ to }  h_i^{d_i}\}$.
Then 
\begin{itemize}
	\item[i)] $Y$ has $$n=\sum_{\substack{1\leq i,j\leq r\\ m_i,\, m_j\mbox{ even}}}
	 \frac{2|G^0|}{m_im_j e_i} \cdot \delta_{ij}$$	
	 singularities of type $A_1$.	 
	\item[ii)] $X$ has $$t= \sum_{\substack{1\leq i\leq r\\ m_i \mbox{ even}}} \frac{|F_i|}{m_i e_i}$$
	 singularities of type $A_3$ and $(n-t)/2$ singularities of type $A_1$.
	\item[iii)] Let $z\in Y$ be a ramification point for $\pi$, 
							then $Q(z)\in \{(p_i,p_i)\mid 1\leq i \leq r\}\subset\bbP^1\times \bbP^1$.
\end{itemize}
\end{prop}

\begin{proof}
i) Let $\xi \in G^0\,, \xi\neq 1$ and assume that $\xi(x,y)=(\xi x, \varphi(\xi) y)=(x,y)$, that is
$$\xi(g K_i, g' K_j)=(g K_i, g' K_j)\Longleftrightarrow
\left\{
\begin{array}{l}
\xi\in g K_i g^{-1}\\
\varphi(\xi) \in g' K_j g'^{-1}
\end{array}
\right .
$$

\noindent So $\tr{\xi}= I= g K_i g^{-1} \cap \varphi^{-1}(g' K_j g'^{-1})\cong \bbZ_2$:
$$\xi=g h_i^{d_i}g^{-1}= \varphi^{-1}(g')\varphi^{-1}(h_j^{d_j})\varphi^{-1}(g'^{-1})\,,$$
for $m_i,\,m_j$ even.

Each element of the form $g h_i^{d_i}g^{-1}$ could stabilize more than one point of $C$, how many?
Let $Z=\{f \in G^0 \mid g h_i^{d_i}g^{-1}=f h_i^{d_i}f^{-1} \}$, since $|Z|=|Z(h_i^{d_i})|=\frac{|G^0|}{e_i}$, 
then $gh_i^{d_i}g^{-1}$ stabilizes $|\{gK_i \mid g \in Z \}|= \frac{|G^0|}{e_i |K_i|}$ points.\\
Each conjugate to $h_i^{d_i}$ stabilizes exactly $\frac{|G^0|}{m_i e_i}$ points in $c^{-1}(p_i)$
and each conjugate to $\varphi^{-1}(h_j^{d_j})$ stabilizes exactly $\frac{|G^0|}{m_j e_i}$ points in $c^{-1}(p_j)$ ($e_i=e_j)$.
Hence if $h_i^{d_i}$ and $\varphi^{-1}(h_j^{d_j})$ are conjugated, there are exactly
$$e_i\cdot \frac{|G^0|}{m_i e_i}\cdot \frac{|G^0|}{m_j e_i}=\frac{|G^0|^2}{m_i m_j e_i}$$
points in $(Q \circ \sigma)^{-1}(p_i,p_j)$ with non-trivial stabilizer.
The orbit of a point $(x,y)$ \sts $\sigma(x,y)$ is a node, has cardinality $\frac{|G^0|}{2}$, hence there are 
$$\frac{|G^0|^2}{m_i m_j e_i}\cdot \frac{2}{|G^0|}= \frac{2|G^0|}{ m_i m_j e_i}$$
nodes on $Y$ over $(p_i,p_j)$.

\

iii) Let $z=\sigma(x,y) \in Y$ be a ramification point for $\pi$ then  $\sigma(x,y)=\sigma(\tau'(x,y))$, 
for some $\tau'\in G$.
If $(Q \circ \sigma)(x,y)=(c(x),c(y))=(p_i,p_j)$ then
$ (Q\circ \sigma)(\tau'(x,y))=(c(y),c(\tau x))=(p_j,p_i)$. Hence $p_i=p_j \in \bbP^1$.

\

ii) We assume $i=1$ and we forget the subscripts.
 
 \noindent  Let $(x,y)\in (Q \circ \sigma)^{-1}(p,p)$  be a ramification point 
for $\pi$ with stabilizer (in $G^0$) $I:= g_1 K g_1^{-1}\cap \varphi^{-1}(g_2 K g_2 ^{-1})\cong \bbZ_2$ then
$$
(g_1 K, g_2 K)=  \tau' \eta (g_1K, g_2K) 
\Leftrightarrow
\left\{
\begin{array}{l}
g_1^{-1}\varphi(\eta) g_2 \in K\\
\varphi(\eta)\tau \eta g_1K =g_1 K \Leftrightarrow (\tau'\eta)^2 \in  g_1K g_1^{-1}
\end{array}
\right.
$$
\noindent Let us fix $\eta \in G^0$ \sts the second condition is fulfilled: 
$o(\tau'\eta)=4$ and $(\tau'\eta)^2=g_1 h^{d} g_1^{-1}$.
Fixed $g_1$, all $|K|$ choices for $g_2$ 
 give the same same point $y\in c^{-1}(p)$.

It can happen that $(\tau'\eta)^2=g_1 h^{d} g_1^{-1}={g'_1} h^{d} {g'_1}^{-1}$ but
$g_1K\neq {g'_1}K$, how many times?
Let $Z=\{u \in G^0 \mid u h^{d} u^{-1}= (\tau'\eta)^2 \}$, since $|Z|=|Z(h^{d})|=\dfrac{|G^0|}{e}$ then 
$\tau' \eta$ stabilizes $|\{uK \mid u\in Z\}|= \frac{|G^0|}{e|K|}=\frac{|G^0|}{e\cdot m}$ points over $(p,p)$.
Suppose that $\tau'\eta(g_1K, g_2 K)=\tau'\xi(g_1 K, g_2 K)$, since $(\tau'\eta)^2$ and $(\tau'\xi)^2$ have both order 2
and $I \ni (\tau'\eta)^2=(\tau'\xi)^2\neq 1$ we conclude that either $\tau'\eta=\tau'\xi$ or $\tau'\eta=(\tau'\xi)^{-1}$.
Hence there are $\dfrac{|F|}{2}\cdot \dfrac{|G^0|}{m e}$ points over $(p,p)$ stabilized by element of the form $\tau'\eta$.
The points in the same orbit for $G^0$ are sent to the same point of $Y$,
so there are 
$$\frac{|F|}{2}\cdot\frac{|G^0|}{m e}\cdot \frac{2}{|G^0|}= \frac{|F|}{m\, e}$$
ramification points for $\pi$ over $(p,p)$.
\end{proof}

\section{The fundamental group}\label{Pi1}

In this section we show how to compute the fundamental group of a regular mixed \qe surface.
 Let $X$ be a mixed \qe quotient determined by $(h_1,\ldots,h_r)$ and $G$ and let 
 $\psi\colon \bbT(m_1,\ldots ,m_r)\rightarrow G^0$ be the appropriate orbifold homomorphism.
The kernel of $\psi$ is isomorphic to the fundamental group $\pi_1(C)$, and 
the action of $\pi_1(C)$ on the universal cover $\D$ of $C$ extends to a discontinuous action of 
$\bbT$. Let $u\colon \D\rightarrow C$ be the covering map, it is $\psi$-equivariant
and $C/G^0\cong \D/\bbT$.

Fix $\tau'\in G\setminus G^0$; let $\tau=\tau'^2\in G^0$ and let 
$\varphi\in \Aut(G^0)$ defined by $\varphi(h):=\tau'h \tau'^{-1}$. 
Let $\bbH:=\{(t_1,t_2)\in \bbT \times \bbT \mid \psi(t_1)=\varphi^{-1}(\psi(t_2))\}\hookrightarrow\Aut(\D\times\D)$;
$\psi$ is surjective and $\varphi(\tau)=\tau$, hence there exists $t\in \bbT$ \sts
$\tilde\tau:=(t,t)\in \bbH$.
We define 
\begin{eqnarray*}
\tilde \tau':\D\times \D& \longrightarrow& \D\times \D\\
(x,y)&\longmapsto& (y,t\cdot x)
\end{eqnarray*}
it is an element of $\Aut(\D\times\D)$  and $(\tilde\tau')^2=\T$;
 we further define $\tilde\varphi \colon \bbH\rightarrow \bbH $ as the conjugation by $\T'$:
$\tilde\varphi(t_1,t_2)=(t_2,t\cdot t_1 \cdot t^{-1})$.

Let $\bbH=\tr{gen(\bbH) \mid rel(\bbH)}$ be a presentation of $\bbH$, and let
$REL:=\{\tilde\varphi(h)\T'h^{-1}\T'^{-1} \mid h\in gen(\bbH)\}$.
We define $\bbG$ as follows:
$$\bbG:=\tr{gen(\bbH), \T' \mid rel(\bbH), (\T')^2\T^{-1}, REL}\,.$$

\begin{thm}\label{ThmPi1}
Let $S\rightarrow X=(C\times C)/G$ be a regular mixed \qe surface. Then 
$$\pi_1(S)\cong\pi_1\bigg(\frac{C\times C}{G}\bigg)\cong \frac{\bbG}{\Tors(\bbG)}\,.$$
\end{thm}

\noindent We recall that the minimal resolution $S\rightarrow X$ 
of the singularities of $X$ replace each singular point by a tree of smooth 
rational curves, hence, by van Kampen's theorem, $\pi_1(S)=\pi_1(X)$.

\noindent To prove the second part of the theorem we need some lemmas.

\

$\bbH$ is an index 2 subgroup of $\bbG$ and 
we define a left action of $\bbG$ on $\D\times \D$, in the following way:
\begin{equation}
\begin{split}
(h_1,h_2)\cdot (x,y) &=(h_1\cdot x, h_2 \cdot y)\\
\T'(h_1,h_2) \cdot (x,y) &= (h_2\cdot y, (t\cdot h_1)\cdot x)
\end{split}
\qquad \mbox{ for } (h_1,h_2)\in \bbH\,.
\end{equation}

\noindent
 We define the homomorphism $\vartheta\colon\bbG\rightarrow G$:
\begin{equation*}
\begin{split}
\vartheta(h_1,h_2)&=\psi(h_1)=\varphi^{-1}\psi(h_2)\\
\vartheta(\T'(h_1,h_2))&= \tau'\psi(h_1)=\tau'\varphi^{-1}(\psi(h_2))
\end{split}
\qquad \mbox{ for } (h_1,h_2)\in \bbH\,.
\end{equation*}
Let $\UU:=(u,u)\colon \D\times \D\rightarrow C\times C$, it is $\vartheta$-equivariant 
and so
$$\frac{\D\times \D}{\bbG} \cong \frac{C\times C}{G}\,,$$
\noindent moreover, we have the following short exact sequence:
\begin{equation*}
1\longrightarrow \pi_1(C\times C)\longrightarrow \bbG\stackrel{\vartheta}{\longrightarrow} G\longrightarrow 1\,.
\end{equation*}

\begin{rem}\label{iso}
	The $\pi_1(C\times C)$-action on $\D\times\D$ is free, so $\pi_1(C\times C)\cap \Stab(x)=\{1\}$, 
	this gives that the restriction of $\vartheta$ to the stabilizer of a point  $x\in \D\times \D$ is 
	an isomorphism onto $\Stab_G(\UU(x))$.
\end{rem}

\begin{lem}\label{disco}
The $\bbG$-action on $\D\times \D$ is discontinuous, that is:
\begin{enumerate}
	\item[(i)] the stabilizer of each point is finite;
	\item[(ii)] each point of $\D\times \D$ has a neighbourhood $U$
	\sts any element of $\bbG$ not in the stabilizer of $x$ maps $U$ outside itself.
\end{enumerate}
\end{lem}

\begin{proof}
(i) By Remark \ref{iso}, the restriction of $\vartheta$ to the stabilizer of $x$ is injective,
	and so $\Stab(x)$ is finite since $G$ is finite.
	
	(ii) Let $x\in \D\times\D$ and let $y:=\UU(x)\in C\times C$,
	since $G$ is finite and $C\times C$ is Hausdorff
	there exists a neighbourhood $U'$ of $y$ \sts for any element of $g\in G$ not in the stabilizer of $y$, $g(U')\cap U' = \emptyset$.
	Let $V'$ be the connected component of $\UU^{-1}(U')$ that contains $x$.
	There exists a connected neighbourhood $V\subseteq V'$ of $x$
	which is mapped isomorphically by $\UU$ onto its image and $\UU(V)=:U\subseteq U'$ is $\Stab(y)$-invariant and 
	$V$ is $\Stab(x)$-invariant.
	Let $\ov{g}\in \bbG\setminus\Stab_\bbG(x)$, we claim that  $\ov{g}(V)\cap V=\emptyset$:
	$$
 \UU(\ov{g}(V)\cap V)\subseteq \UU(\ov{g}(V))\cap \UU(V)=\vartheta(\ov{g})U\cap U \,,
$$
	hence either $\ov{g}(V)\cap V= \emptyset$ or $\vartheta(\ov{g})\in \Stab(y)$. In the latter case, by Remark \ref{iso}, there exists a unique
	$\ov{g}'\in \Stab(x)$ \sts $\vartheta(\ov{g}')=\vartheta(\ov{g})$, so $\ov{g}= k \ov{g}'$ with
	$k\in \pi_1(C\times C)\setminus\{1\}$ and we get:
	 $$
	 \ov{g}(V)\cap V = k\ov{g}'(V)\cap V= k (V)\cap V=\emptyset\,.
	$$
\end{proof}

\begin{defi}
Let $H$ be a group, its \no{torsion subgroup} $\Tors(H)$ is the normal subgroup generated by
all elements of finite order in $H$.
\end{defi}

\begin{lem}\label{tors}
The normal subgroup $\bbG'$ of $\bbG$ generated by the elements which have non-empty 
fixed-point set is exactly $\Tors(\bbG)$.
\end{lem}

\begin{proof}
To prove our claim we show that each element $g\in \bbG$ of finite order has non-empty 
fixed-point set, and vice versa. We distinguish two cases: 
\begin{itemize}
	\item[(i)] Let $g=(h_1,h_2)$ be an element of $\bbH\subset \bbG$ that fixes a point 
	$(x,y) \in \D\times \D$:
	$$(h_1,h_2)(x,y)=(x,y) \Longleftrightarrow \left\{\begin{array}{c}
	h_1=\alpha c_i^{m_i}\alpha^{-1}\\
	h_2=\beta c_j^{m_j}\beta^{-1}
	\end{array}\right.
	\Longleftrightarrow
	(h_1,h_2) \mbox{ has finite order;}
	$$
	the first equivalence follows by the proof of the Riemann existence theorem, while for the second see \cite[Theorem 10.3.2]{Bear83}.
	
	\item[(ii)] Let $g=\T' (h_1,h_2)\in \bbG\setminus\bbH$.
	If $g$ fixes a point $(x,y) \in \D\times \D$, also $g^2\in\bbH$ fixes the point, by (i) it has
	finite order, then $g$ has finite order.
	Conversely, if $g$ has finite order, $g^2(x,y)=(x,y)$ for some $(x,y) \in \D\times \D$
	since $g^2\in \bbH$ has finite order and
 $g(x,(h_2^{-1}) x)=(x,(h_2^{-1}) x)$.		
\end{itemize}
\end{proof}

\begin{proof}[Proof of Theorem \ref{ThmPi1}]

Because of Lemma \ref{disco}, the main theorem in \cite{Arm68} applies and we get:
$$\pi_1\bigg(\frac{C\times C}{G}\bigg)\cong\pi_1\bigg(\frac{\D\times \D}{\bbG}\bigg)\cong \frac{\bbG}{\bbG'}$$
where $\bbG'$ is the normal subgroup of $\bbG$ generated by the elements which have non-empty 
	fixed-point set, which is exactly $\Tors(\bbG)$ by Lemma \ref{tors}:
	$$\pi_1\bigg(\frac{C\times C}{G}\bigg)\cong \frac{\bbG}{\Tors(\bbG)}\,.$$
\end{proof}

To build a MAGMA script that calculates the 
fundamental group, we have to find a finite set of generators of $\Tors(\bbG)$.

\begin{prop}\label{fingen}
Let $X$ be a mixed \qe quotient determined by $(h_1,\ldots,h_r)$ and $G$ and let 
 $\psi\colon \bbT(m_1,\ldots ,m_r)\rightarrow G^0$ be the corresponding appropriate orbifold homomorphism.
 Fix $\tau'\in G\setminus G^0$; let $\tau=\tau'^2\in G^0$ and let 
$\varphi\in \Aut(G^0)$ defined by $\varphi(h):=\tau'h \tau'^{-1}$. 
 Then 
$\Tors(\bbG)$ is normally generated by the finite set $T_1 \cup T_2$ constructed as follows:

\begin{itemize}

\item $T_1\subset\bbH$: for every $i,j \in\{1,\ldots, r\}$, $1\leq \alpha \leq m_i-1$ and $1\leq \beta \leq m_j-1$,
if $h_i^\alpha$ is conjugated to $\varphi^{-1}(h_j^\beta)$, then we choose
an element $v\in G^0$ such that $v\,h_i^\alpha\, v^{-1}=\varphi^{-1}(h_j^\beta)$.
Then for every element $d$ in the finite group $Z(h_i^\alpha)$ we choose
an element $w\in \psi^{-1}(v\cdot d)$  and we include $(w\,c_i^\alpha\, w^{-1}, c_j^\beta)$  in $T_1$.

\item $T_2\subset\bbG\setminus \bbH$:
for every $i, \in\{1,\ldots, r\}$, $1\leq \alpha \leq m_i-1$
and $\eta\in G^0$, if $(\tau' \eta)^2$ is conjugated to $h_i^\alpha$,
then we choose an element $v\in G^0$ such that $v\,h_i^\alpha\, v^{-1}=(\tau' \eta )^2$ and 
we choose  $g_1\in \psi^{-1}(\eta)$ and $g_2\in \psi^{-1}(\varphi(\eta))$.
Then for every element $d$ in the finite group $Z(h_i^\alpha)$ we choose
an element $w\in \psi^{-1}(v\cdot d)$,
and we include $\T'(g_1,kg_2)$ in $T_2$, where $k:=(g_2\,t\,g_1)^{-1}w c_i^\alpha w^{-1}$.
\end{itemize}
\end{prop}

\begin{proof}
By \cite[Lemma 4.9]{BCGP08}, $T_1$ normally generates $\Tors(\bbH)$ which is also
the set of the elements of $\bbH$ that stabilize some points in $\D\times \D$.

Let $\eta\in G^0$ such that $\tau'\eta(x,y)=(x,y)$ for some $(x,y)\in C\times C$, i.e.

$$
\tau'\eta(x,y)=(x,y)\quad \Longleftrightarrow \quad
\left\{
\begin{array}{l}
x=\varphi(\eta) \,y\\
y=\tau \eta\, x
\end{array}
\right.
\quad\Longleftrightarrow\quad
\left\{
\begin{array}{l}
x=(\tau' \eta)^2x\\
y=\tau \eta  \,x
\end{array}
\right.
$$
So $\tau'\eta $ stabilizes some points in $C	\times C$ if and only if $(\tau' \eta)^2$ is 
conjugated to $h_i^\alpha$ for some $1\leq i\leq r $ and $1\leq \alpha \leq m_i-1$.

\noindent Fix $g_1\in \psi^{-1} (\eta)$ and $g_2 \in \psi^{-1}(\varphi(\eta))$,  the preimages of $\tau'\eta$  are of the form
$\T'(g_1k_1,g_2k_2)$, where $k_1,\,k_2\in \ker \psi$,
but up to conjugation with $(k_1,1)\in \bbH$, we can assume that
they are of the form $\T'(g_1,kg_2)$ with $k\in \ker \psi$.

\noindent Let $s:=g_2 t g_1 \in \bbT$: $$\psi(s)=(\tau'\eta)^2=v h_i^\alpha v^{-1}$$ for some $v \in G^0$.
For any $d \in Z(h_i^\alpha)$, 
let $w$ be a preimage of $v\cdot d$ via $\psi$, so $s=w c_i^\alpha w^{-1}k'$ where $k'\in \ker \psi$.
We define $$k:=(k')^{-1}=s^{-1} w c_i^\alpha w^{-1}\,,$$
 hence $ks$ is conjugated to $c_i^\alpha$ and so it stabilizes some point $x_0\in\D$ and 
 $\T'(g_1,kg_2)$  stabilizes $(x_0, (kg_2)^{-1} x_0)\in\D\times\D$,
moreover $\UU(x_0, (kg_2)^{-1} x_0)=(x,y)$.
We include $\T'(g_1,kg_2)$ in $T_2$.

We are left with showing that every element
in $\bbG\setminus \bbH$ that stabilizes some points in $\D\times \D$ belongs to the subgroup 
normally generated by $T_1\cup T_2$.

Let $\T'(h_1,h_2)\in \bbG$ be an element that stabilizes a point $(x_1,y_1)\in \D \times \D$.
There exists $g\in T_2$ and $(x_0,y_0)\in \UU^{-1} (\UU(x_1,y_1))$ such that $g(x_0,y_0)=(x_0,y_0)$.
 By construction, there exists $g'\in \bbG$ such that
$g'(x_0,y_0)=(x_1,y_1)$, hence $g'gg'^{-1}(x_1,y_1)=(x_1,y_1)$. 
By remark \ref{iso}, there exists $h\in \Tors(\bbH)$ such that $\T'(h_1,h_2)= h g'gg'^{-1}$.
\end{proof}

\section{The classification of the surfaces}\label{class}

In this section we give a complete classification of the regular 
mixed \qe surface $S$ with $p_g(S)=0$ occurring as the minimal resolution of the singularities of a  mixed \qe quotient
$X := (C\times C)/G$  such that 	 $(C\times C)/G^0$ has only nodes as singularities.

We let the computer make a systematic search of the surfaces that satisfy the above assumptions.
As output we get the following theorem:

\begin{thm}
Let $S$ be a minimal regular surface with $p_g(S)=0$ whose canonical
model is the mixed \qe quotient $X=(C\times C)/G$ \sts $(C\times C)/G^0$ has at most nodes as singularities.
Then $S$ is of general type and belongs to one of the 13 families collected in Table \ref{tab1}.
\end{thm}
 
\begin{rem} \label{missing}
In \cite{BCG08} the authors have considered the case in which the quotient map $C\times C\rightarrow X$ 
is \'etale, \ie $G$ acts freely on $C\times C$. In this case $X$ is smooth and $K_S^2=8$. 
Running the MAGMA script in this special case we noted that they missed
 a family of surfaces, that is tagged by {A.5.1}.
\end{rem}

\begin{rem}\label{qhpp}
For the canonical model of a surface of general type with $p_g=0$
is automatic that $b_0=b_4=1$ and $b_1=b_3=0$, hence those with $b_2=1$
are $\bbQ$-homology projective planes, see \cite{HK} and \cite{Keum10}.
The surfaces labelled by {A.3.1} and {A.3.2} are new examples of $\bbQ$-homology projective planes.
\end{rem}

\begin{rem}\label{z4}
 We point out that the surface {A.3.1} is a numerical Campedelli surface ($K_S^2=2$) with topological fundamental group
	(and therefore algebraic fundamental group) $\bbZ_4$. We have discussed the importance of this surface in the introduction.
\end{rem}

	\begin{rem}\label{news}
	We have constructed  2 new topological types of surfaces of general type with $p_g=0$.
	These surfaces are tagged by {A.4.1} and {A.4.3}.
\end{rem}
\begin{rem}\label{K-N}
 The surface tagged by {A.4.2} has $K_S^2=4$ and the same fundamental group of a 
	 Keum-Naie surface (see \cite{Naie} and \cite{BC09}).
We expect that this surface belongs to the family studied in \cite{BC09} but we have not proved it.
\end{rem}
\begin{rem}
There has been a growing interest for surfaces of general type with $p_g=0$ having an involution,
see \cite{CCMM}, \cite{CMMP}, \cite{Rit10} and \cite{YY10}.
The \q{intermediate} surface $Y=(C\times C)/G^0$ has an involution given by $\sigma\colon Y\rightarrow X$;
it has $q=0$ and $K^2_Y=2K^2_S$, while
 $p_g=0$ in the cases {A.1.1}, {A.3.1} and {A.3.2}, and $p_g=1$ in the others.

Let $S\rightarrow X$ be the Godeaux surface  ($K^2_S=1$) tagged by {A.1.1}.
 The surface $Y$ has 6 nodes and $K^2_Y=2$, moreover its desingularization $T$
inherits an involution $\nu$ from the involution acting on $Y$ and has $K^2_T$, hence we have a Campedelli surface with an involution.
By construction, the involution fixes 4 points on $T$, by \cite[Proposition 2.3]{CMMP} in this case the involution is not composed
with the bicanonical map $\varphi \colon T\rightarrow\bbP^2$. By construction $S$ is also the desingularization of $T/\tr{\nu}$, 
this means that $S$ is an example of the case (i) of \cite[Proposition 4.3]{CMMP}.

In the cases {A.3.1} and {A.3.2}, $Y$ is a surface with $K_Y^2=4$, $p_g=0$ and $4$ nodes. 
These surfaces are the quotient models of two product-quotient surfaces constructed in
\cite{BCGP08}.
\end{rem}

\subsection{The script}

Using the results of the previous sections we implement a MAGMA script to find all the
surfaces satisfying our assumptions.
The algorithm follows closely the algorithms in \cite{BCGP08} and \cite{BP10}. 
We have extended them to the mixed case and we have improved the computational complexity.
We explain the strategy of the program and  the most important scripts;
a commented version of the full code can be downloaded from:
$$\mbox{\url{http://www.science.unitn.it/~frapporti/papers/scriptmix1.magma}}\,$$

First of all we fix a value of $K_S^2\in\{1,\ldots, 8\}$.
\begin{enumerate}
	\item[Step 1:] the script \textbf{Sings} list all the possible baskets of singularities for 
	$K_S^2$, accordingly with Lemma \ref{K2value} there are only finitely many.
	
	\item[Step 2:] by Lemma \ref{inequalities}, once we fix $K_S^2$
	there are finitely many possible signatures.
	\textbf{ListOfTypes} computes them. The input is $K_S^2$, so this script before computes $Sings(K_S^2)$ and 
	returns a list of pairs: the first entry is a possible basket and the second is the list with all the possible signatures.	
	
	\item[Step 3:] if we know the signature, by Proposition \ref{orderG}, we can compute the order of $G^0$.
	\textbf{ListGroups}, whose input is $K_S^2$, searches, for every element in the output of $ListOfTypes$, if among 
	the group of the right order there are groups having at least one set of spherical generators of the 
	prescribed type. Then it checks if these groups have a pair of set of spherical generators 
	that give the prescribed basket of singularities on $Y=(C\times C)/G^0$.	
	Once it finds a group $G^0$ with the right properties, it searches among all the groups of order $2|G^0|$
	 the  unsplit extensions of $G^0$.
		
\noindent	For each positive answer to these questions $ListGroups$ stores the triple $(basket, type, (group, j))$,  
where group is the group $G$ and $j$ identifies	the group $G^0$ as subgroup of $G$.	

	The script has some conditional instructions:
	
\begin{itemize}

	\item if one of the signatures is $(2,3,7)$, then $G^0$, being a quotient of $\bbT(2, 3, 7)$,
is perfect. MAGMA knows all perfect groups of order $\leq 50000$, and then
$ListGroups$ checks first if there are perfect group of the right order: if not,
this case can not occur.

 \item If the order is a number as e.g., 576, where there are too many
isomorphism classes of groups, then $ListGroups$ makes the controls in a different way, i.e it
use the MAGMA function \q{SmallGroupProcess} and not the function \q{SmallGroups} to find the 
groups of order $|G^0|$.

\item If the order of $G^0$ is in $\{1001,\ldots, 2000\}$, since MAGMA does not have a list
			of the groups of order bigger than 2000, 	
			we can not check if there exist groups that extends it in a non split way.
			We collect these cases in a list, second output of the script.

\item  if the expected order of the group $G^0$ is 1024 or bigger than 2000,
since MAGMA does not have a list of the finite groups of these orders,
 then $ListGroups$ just stores these cases in a list, third output of the script.
	\end{itemize}
	\item[Step 4:] \textbf{ExistingSurfaces} takes the output of $ListGroups(K_S^2)$ and throws away all triples that do not give a 
	surfaces with the expected singularities.
	
	\item[Step 5:] each triple in the output of $ExistingSurfaces(K_S^2)$ gives many surfaces, one for each spherical systems of generators.
	Two different spherical systems of generators can give deformation equivalent surfaces (see Remark \ref{cl}),
		the script \textbf{FindSurfaces} produces one representative for each equivalence class.
	
	\item[Step 6:] \textbf{Pi1} computes the fundamental group of the surfaces.
	\end{enumerate}

\begin{rem}\label{abel}
The principal computational improvement in our script is in the first part of $ListGroups$, 
in particular in the search of which groups have at least a set of spherical generators of the 
	prescribed type: if $G^0$ has a set of spherical generators of type $(m_1,\ldots,m_r)$,
	then there exists an appropriate orbifold homomorphism $\psi \colon \bbT(m_1,\ldots,m_r)\rightarrow G^0$.
	The map $\psi$ induces a surjective morphism $\ov{\psi}\colon \bbT^{ab}\rightarrow {G^0}^{ab}$ between the
	abelianizations, hence ${G^0}^{ab}$ is isomorphic to a quotient of $\bbT^{ab}$. 
	Our script checks first  (by the script {\bf Test}) which groups have abelianization 
	isomorphic to a quotient of the suitable $\bbT^{ab}$ and 
	only for the groups that pass this test it search for a set of spherical generators of the right type.
	\end{rem}		
 In the table below we compare the execution times of a script with $Test$ and one without $Test$.
 We note that using $Test$ the execution time decreases up to 20 times.	
	
\begin{center}
\begin{table}[!ht]
\begin{tabular}{|c|c|c|c|c|c|}
\hline
 $K^2$& 1&2&3&4&5\\
 \hline
\textbf{Test}  &12.59s& 983.54s &109.31s &296.16s & 31.95s \\
 \hline
\textbf{NO Test} &111.96s&4307.94s &2281.98s & 4740.24s&237.96s  \\
 \hline
\end{tabular}
\caption{Some execution times}
\end{table}
\end{center}

\subsection{Skipped Cases}

\

The script \textbf{ListGroup} returns 3 output: the first is processed by the other scripts that possibly return some surfaces.
The other outputs are cases which we have to study separately. 
  We will show that all these cases do not occur.

For all the values $1\leq K^2_S\leq 8$, we have that the second output is empty, while the cases stored in
the third output are collected in Table \ref{tabSkip}.

\begin{center}
\begin{table}[!ht]
\begin{tabular}{|c|c|c|c|}
\hline
$K^2_S$ & $\Sing X$ & type & $|G^0|$\\
\hline
4 & $4 \times A_1$ &   2, 3, 8  & 2304 \\
\hline
5& $3\times A_1$ &   2, 3, 8   &2880\\
\hline
6 & $2 \times A_1$ &   2, 4, 5   & 2400\\
6 & $2 \times A_1$ &   2, 3, 8   & 3456\\
\hline
\multicolumn3c{}\\
\end{tabular}
\qquad
\begin{tabular}{|c|c|c|c|}
\hline
7 & $1\times A_1$ &   2, 3, 9   & 2268  \\
7 & $1\times A_1$ &   2, 4, 5   & 2800  \\
7 & $1\times A_1$ &   2, 3, 8   & 4032  \\
\hline
8 & $\emptyset$ &   2, 3, 9   & 2592\\
8 & $\emptyset$ &   2, 4, 5   & 3200\\
8 & $\emptyset$ &   2, 3, 8   & 4608\\
\hline
\end{tabular}
\caption{The  cases skipped by $ListGroups$}
\label{tabSkip}
\end{table}
\end{center}

Here we show only the case $K^2=7$, $|G|=4032$ that contains all the
tools used to prove that these cases do not occur.
The proofs for the
other cases are completely analogous and can be found at
the following link: $$\mbox{\url{http://www.science.unitn.it/~frapporti/papers/skipped1.pdf}}\,$$

\begin{lem}
No group of order 4032 
has a pair of spherical system of generators of type $(2,3,8)$ which gives the expected singularities, i.e. 2 nodes on $Y$.
\end{lem}

\begin{proof}
Assume that $G^0$ is a group of order 4032 with a spherical system of generators of type $(2,3,8)$: $(a,b,c)$.
Since $\bbT(2,3,8)^{\mbox{ab}}\cong \bbZ_2$ and  since there are no perfect groups of order 4032:
\begin{verbatim}
> NumberOfGroups(PerfectGroupDatabase(),4032);
0
> \end{verbatim}
the commutator subgroup ${G^0}'=[G^0,G^0]$ of $G^0$ has order 2016.
Using the Reidemeister-Schreier method (see \cite[Section 2.3]{MKS}) or an easy  MAGMA function,
one can easily  show that  $[\bbT(2,3,8),\bbT(2,3,8)]\cong\bbT(3,3,4)$, moreover 
$(d,e,f):=(aba, b, c^2)$  is a spherical set of generators of type $(3,3,4)$ for ${G^0}'$.

Since $\bbT(3,3,4)^{\mbox{ab}}\cong \bbZ_3$ and  since there are no perfect groups of order 2016,
 the commutator subgroup ${G^0}''=[{G^0}',{G^0}']$ of ${G^0}'$ has order 672,  it is a quotient of 
$[\bbT(3,3,4),\bbT(3,3,4)]\cong\bbT(4,4,4)$ and $(efe^{-1}, e^2fe^{-2}, f)$ is a spherical set of generators of type $(4,4,4)$ for ${G^0}''$.
The following MAGMA computation
\begin{verbatim}
> Test([4,4,4], 672);
{ 1046, 1255 }
>
\end{verbatim}
shows that only the groups $G(672,v)$ with $v\in\{1046,1255\}$ 
have a spherical system of generators of type $(4,4,4)$.
By assumption, ${G^0}'$ has a spherical system of generators of type $(3,3,4)$, then ${G^0}' \cong {G^0}''\rtimes \bbZ_3$.
The next claim is a standard result about semidirect product.

\begin{claim}
Let $L$ be a finite group and let $K$ be a cyclic group of order $p$.
Let $\varphi_1 ,\varphi_2:K \rightarrow \Aut(L)$ \sts $ \varphi_1(K)$ and $\varphi_2(K)$ are conjugated.
Then $L\rtimes_{\varphi_1}K \cong L\rtimes_{\varphi_2}K$.
\end{claim}

In order to construct the group ${G^0}'$, we have only to consider the conjugacy classes of 
elements of order 3 in $\Aut({G^0}'')$ and $\Id(\Aut({G^0}''))$.

The following MAGMA script shows that ${G^0}''=G(672,1046)$
 has only one conjugacy class of automorphisms of order 3; 
hence, up to isomorphisms, there are at most two ${G^0}'\cong{G^0}''\rtimes \bbZ_3$.
The script shows also that these two extensions ${G^0}''\rtimes \bbZ_3$ do
 not have a spherical system of generators of type $(3,3,4)$, hence these cases do not occur.
\begin{verbatim}
> H2:=SmallGroup(672,1046);
> R2:=AutConjugCl(H2,3);
56
1
> C3:=CyclicGroup(3);
> Aut2:=AutomorphismGroup(H2);
> R2[2]:=Id(Aut2);
> f:=[]; for i in [1..2] do f[i]:=hom<C3->Aut2|R2[i]>;end for;
> h1:=[]; for i in [1..2] do h1[i]:=SemidirectProduct(H2,C3,f[i]);
for> i, ExSphGens(h1[i],[3,3,4]); end for;
1 false
2 false
>
\end{verbatim}

The following MAGMA script shows that ${G^0}''=G(672,1255)$ has,
up to isomorphisms, only one  extension ${G^0}' \cong {G^0}''\rtimes \bbZ_3$
with a spherical system of generators of type $(3,3,4)$.

\begin{verbatim}
> H2:=SmallGroup(672,1255);
> R2:=AutConjugCl(H2,3);
170
3
> C3:=CyclicGroup(3);
> Aut2:=AutomorphismGroup(H2);
> R2[4]:=Id(Aut2);
> f:=[]; for i in [1..4] do f[i]:=hom<C3->Aut2|R2[i]>;end for;
> h1:=[]; for i in [1..4] do h1[i]:=SemidirectProduct(H2,C3,f[i]);
for> i, ExSphGens(h1[i],[3,3,4]); end for;
1 true
2 false
3 true
4 false
> IsIsomorphic(h1[1],h1[3]);
true  Homomorphism of ...
>H1:=h1[1];
\end{verbatim}

By assumption, $G^0$ has a spherical system of generators of type $(2,3,8)$, then $G^0\cong{G^0}'\rtimes \bbZ_2$.
The following MAGMA script (that continues the previous one) shows that ${G^0}'=h1[1]$ has,
up to isomorphisms, only one  extension ${G^0}'\rtimes \bbZ_2$
with a spherical system of generators of type $(2,3,8)$.
\begin{verbatim}
> R1:=AutConjugCl(H1,2);
499
8
>C2:=CyclicGroup(2);
> Aut1:=AutomorphismGroup(H1);
> R1[9]:=Id(Aut1);
>f:=[]; for i in [1..9] do f[i]:=hom<C2->Aut1|R1[i]>;end for;
> h:=[]; for i in [1..9] do h[i]:=SemidirectProduct(H1,C2,f[i]);
for> i, ExSphGens(h[i],[2,3,8]); end for;
1 true
2 false
3 false
4 false
5 false
6 false
7 false
8 false
9 false
>
\end{verbatim}

The following MAGMA script shows that for ${G^0}=h[1]$ the singularities test fails,
and so  this case does not occur.
  
\begin{verbatim}
>H:=h[1];
> SingularitiesY([0,1],[2,3,8],H);
false
\end{verbatim}
\vspace{-0.5cm}
\end{proof}


\noindent
\textbf{Acknowledgments} The author thanks a lot his advisor R. Pignatelli  
 for a lot of useful discussions and suggestions.
 The author also thanks D. Hwang for pointing out a miscomputation in a previous version
 of the paper.


\

\

\noindent Davide Frapporti\\
Dipartimento di Matematica della Universit\`a di Trento;\\
Via Sommarive 14; I-38123 Trento (TN), Italy\\
e-mail: {\verb* davide.frapporti@gmail.com }
\end{document}